\documentclass[twoside,10pt,leqno]{amsart}
\everymath={\displaystyle}
\usepackage{amsfonts}
\usepackage{amsmath}
\usepackage{amscd}
\usepackage{amssymb}
\usepackage{amsthm}
\usepackage{amsrefs}
\begin{document}
\newtheorem{theorem}{Theorem}
\newtheorem{lemma}{Lemma}
\newtheorem{corollary}{Corollary}
\newtheorem{conjecture}{Conjecture}
\newtheorem{prop}{Proposition}
\numberwithin{equation}{section}
\newcommand{\dif}{\mathrm{d}}
\newcommand{\intz}{\mathbb{Z}}
\newcommand{\ratq}{\mathbb{Q}}
\newcommand{\natn}{\mathbb{N}}
\newcommand{\comc}{\mathbb{C}}
\newcommand{\rear}{\mathbb{R}}
\newcommand{\prip}{\mathbb{P}}
\newcommand{\uph}{\mathbb{H}}

\title{\bf Large Sieve Inequalities for Special Characters to Prime Square Moduli}
\author{Liangyi Zhao}
\maketitle

\begin{abstract}
In this paper, we develop a large sieve type inequality for some special characters whose moduli are squares of primes.  Our result gives non-trivial estimate in certain ranges.
\end{abstract}

\section{Introduction and historical background}

It was in 1941 that Yuri Vladimirovich Linnik \cite{JVL1} originated the idea of large sieve.  It is as follows.  A set of real numbers $\{ x_k \}$ is said to be $\delta$-spaced modulo 1 if $x_j-x_k$ is at least $\delta$ away from any integer, for all $j \neq k$. \newline

The large sieve inequality, which we henceforth refer to as the classical large sieve inequality, is stated as follows.  Different elegant proofs of the theorem can be found in \cite{HD}, \cite{PXG}, \cite{HM2}, \cite{HM}.  The theorem, in the following form, was first introduced by Davenport and Halberstam, \cite{DH1} and \cite{DH2}.

\begin{theorem}
Let $\{ a_n \}$ be an arbitrary sequence of complex numbers, $\{ x_k \}$ be a set of real numbers which is $\delta$-spaced modulo 1, and $M \in \intz$, $N \in \natn$.  Then we have
\begin{equation*}
\sum_k \left| \sum_{n=M+1}^{M+N} a_n e \left( x_k n \right) \right|^2 \ll \left( \delta^{-1} +N \right) \sum_{n=M+1}^{M+N} |a_n|^2,
\end{equation*}
where the implied constant is absolute.
\end{theorem}
Save for the more precise implied constant, the above inequality is the best possible.  This theorem admits corollaries for multiplicative characters.  We derive
\begin{equation} \label{lssmult}
\sum_{q=1}^Q \frac{q}{\varphi(q)} \; \sideset{}{^{\star}}\sum_{\chi \; \bmod{ \; q}} \left| \sum_{n=M+1}^{M+N} a_n \chi (n) \right|^2 \ll (Q^2+N) \sum_{n=M+1}^{M+N} |a_n|^2,
\end{equation}
where here and after, $\sideset{}{^{\star}}\sum$ means that the sum runs over primitive characters modulo the specified modulus only.  Note that the sum is over only the primitive characters is a vital feature of the theorem. \newline

Also, large sieve type inequalities for the complete set of characters to square, and higher power,  moduli was investigated by L. Zhao in \cite{LZ1}. \newline

In this paper, we aim to have a result of the following kind.
\begin{equation} \label{lsspecini}
\sum_{q=1}^Q \frac{q}{\varphi(q)} \; \sideset{}{^{'}}\sum_{\chi \bmod{q^2}} \left| \sum_{n=M+1}^{M+N} a_n \chi (n) \right|^2 \ll \Delta \sum_{n=M+1}^{M+N} |a_n|^2,
\end{equation}
where the $\sideset{}{^{'}}\sum$ runs over some special Dirichlet characters, to be specified in the next section, and as usual $\Delta$ will be in terms of $Q$ and $N$. \newline

In all of our investigations of this paper, we shall restrict our attention to prime square moduli only.  The result can be generalized, but restricting the prime square moduli gives us great convenience in estimates, {\it id est} two prime squares are not co-prime if and only if they are the same. \newline

The idea of studying character sums to prime-power moduli was started in a paper by A. G. Postnikov \cite{AGP}.  He gave formulas on the decomposition of groups of characters of powerful moduli.  P. X. Gallagher \cite{PXG2} also studied characters of this type.  Iwaniec \cite{HI6} expanded the Postnikov-Gallagher idea to composite moduli. \newline

The author also thanks his thesis adviser, Henryk Iwaniec, who first suggested this problem and who, of his advise and support, has been most generous.

\section{Heuristics}

First we note that there are $\varphi(q^2) = q\varphi(q)$ Dirichlet characters modulo $q^2$ which is the same as the order of the group $G=\left( \intz / q^2 \intz \right)^*$.  The group $G$ contains the following subgroup $H = \left\{ x \in G \; | \; x \equiv 1 \pmod{q} \right\}$, which is isomorphic to the additive group $\intz /q \intz$.  The isomorphism is given by
\[ H \longrightarrow \intz /q \intz : x \longrightarrow \frac{x-1}{q}. \]
Any Dirichlet character $\xi$ on $G$ induces an additive character on $H$.  So we have, for any $x \in H$,
\begin{equation} \label{isomorph}
\xi (x) = e \left[ \frac{a(x-1)}{q^2} \right],
\end{equation}
for some $a \pmod{q}$.  If $\xi' \pmod{q^2}$ is another character satisfying \eqref{isomorph} with the same $a \pmod{q}$.  Then $\xi' \xi^{-1}$ is a character on $G$ that is trivial on $H$.  Therefore, $\xi'=\xi \chi$ where $\chi$ is induced by a character on $\left( \intz / q \intz \right)^*$.  Let $G_a$ be the set of characters $\xi \pmod{q^2}$ satisfying \eqref{isomorph}.  $G_a$ has $\varphi(q)$ elements, and that every element of $G_a$ is obtained in a unique way by multiplying a fixed character $\xi \in G_a$ by a character $\chi$ on $\left( \intz / q \intz \right)^*$. \newline

Therefore, the context in which we shall consider the sum in the left-hand side of \eqref{lsspecini} is with $a=1$ in \eqref{isomorph}; {\it Id est}, we shall study the following sum
\begin{equation*}
\sum_{q=1}^Q \frac{q}{\varphi(q)} \sum_{\xi \in G_1} \left| \sum_{n=M+1}^{M+N} a_n \xi (n) \right|^2.
\end{equation*}
The more general case in which one considers characters in $G_a$ is similar to the investigation for $G_1$.  Also, we note that the characters being summed are not necessarily primitive, a feature that the classical large sieve inequality, \eqref{lssmult}, does not possess. \newline

Next, we observe the following.
\begin{equation*}
\frac{q}{\varphi(q)} \sum_{\xi \in G_1} \left| \sum_{n=M+1}^{M+N} a_n \xi (n) \right|^2
= \frac{q}{\varphi(q)} \sum_n \sum_{n'} \overline{a}_n a_{n'} \xi(\overline{n} n') \sum_{\chi} \chi (\overline{n} n'),
\end{equation*}
where $\xi$ is a fixed character in $G_1$ and $\chi$ runs over Dirichlet character of $\left( \intz / q \intz \right)^*$.  The inner-most sum vanishes unless $n \equiv n' \pmod{q}$, in which case it yields $\varphi(q)$ and
\begin{equation*}
\xi (\overline{n}n')= e \left( \frac{\overline{n}n'-1}{q} \right) = e \left( \frac{\overline{n}}{q} \frac{n'-n}{q} \right)
\end{equation*}
by \eqref{isomorph}.  Note that $\frac{n'-n}{q}$ is an integer.  Hence, we have
\begin{equation} \label{inieq}
\frac{q}{\varphi(q)} \; \sideset{}{^{(1)}} \sum_{\xi \bmod{q^2}} \left| \sum_{n=M+1}^{M+N} a_n \xi (n) \right|^2 = q \mathop{\sum_n \sum_{n'}}_{n \equiv n' \bmod{q}} \overline{a}_n a_{n'} e \left( \frac{\overline{n}}{q} \frac{n'-n}{q} \right),
\end{equation}
where henceforth $\sideset{}{^{(1)}} \sum$ denotes sum over characters $\xi \in G_1$.  Trivially estimating the contribution of the above sum gives
\begin{equation} \label{lsspectriv}
\frac{q}{\varphi(q)} \; \sideset{}{^{(1)}} \sum_{\xi \bmod{q^2}} \left| \sum_{n=M+1}^{M+N} a_n \xi (n) \right|^2 = \left[ q + O\left( N \right) \right] \sum_{\gcd(n,q)=1} |a_n|^2.
\end{equation}
Estimating thus gives that
\begin{equation} \label{lsspectriv2}
\sum_{q=1}^Q \frac{q}{\varphi(q)} \; \sideset{}{^{(1)}} \sum_{\xi \bmod{q^2}} \left| \sum_{n=M+1}^{M+N} a_n \xi (n) \right|^2 \ll (Q^2+QN) \sum_n |a_n|^2.
\end{equation}
Although \eqref{lsspectriv} is obtained trivially, it is an asymptotic formula rather than simply an upper bound, a main feature of the classical large sieve inequalities.  The inequality \eqref{lsspectriv2} is of interest when $N \ll Q$.  Moreover, in the light of \eqref{lsspectriv}, we see that any improvement upon the exponent of $Q$ is not possible and any improvement upon this trivial bound has to come from that on $N$.  Toward that end, we shall use the results in the following section. \newline

\section{Preliminaries Lemmas}

First, we shall need the estimate for Kloosterman sums.

\begin{theorem}[Weil] \label{weilbound}
For any $c \in \natn$, $m$ and $n$, we have
\begin{equation*}
\left| \sum_{ad \equiv 1 \bmod{c}} e \left( \frac{ma+nd}{c} \right) \right| \leq [\gcd(m,n,c)]^{\frac{1}{2}} c^{\frac{1}{2}} \tau(c),
\end{equation*}
where $\tau(c)$ is the divisor function.
\end{theorem}
\begin{proof} This is quoted from \cite{HI3} and is deduced from the celebrated Riemann hypothesis for curves over finite fields proved by A. Weil in 1948. \end{proof}

Next we shall also need the following estimate for Ramanujan's Sums.

\begin{lemma} For any $n,q \in \natn$, we have 
\begin{equation} \label{ramsumest}
\left| \sum_{\substack{a \pmod{q} \\ \gcd(a,q)=1}} e \left( \frac{an}{q} \right) \right| \leq \gcd(n,q).
\end{equation}
\end{lemma}

\begin{proof} This estimate is easy, standard, and best possible. \end{proof}

\section{Main Contention}

We now state and prove the theorem of this paper.  Throughout, $q$ runs over prime numbers only.

\begin{theorem} \label{lsspec}
Suppose $Q, N \in \natn$, $M \in \intz$ and $\{ a_n \}$ be a sequence of complex numbers.  We have
\begin{equation} \label{lsspeceq}
\sum_{q=1}^Q \frac{q}{\varphi(q)} \; \sideset{}{^{(1)}} \sum_{\xi \bmod{q^2}} \left| \sum_{n=M+1}^{M+N} a_n \xi (n) \right|^2 \ll Q^{\epsilon} \left( N Q^{\frac{1}{2}} + N^{\frac{1}{4}} Q^2 + N^{\frac{3}{4}} Q^{\frac{11}{8}} \right) \sum_{n=M+1}^{M+N} |a_n|^2,
\end{equation}
with any $\epsilon >0$ and the implied constant depends on $\epsilon$ alone.
\end{theorem}

\begin{proof}
For the proof, it suffices to consider $q \in \left( Q/2,Q\right]$.  The size of th sum resulted from breaking the left-hand side of \eqref{lsspeceq} into dyadic intervals is majorized by $QT(N,Q)$, with
\begin{equation*}
T(N,Q) = \sum_{Q/2<q\leq Q} \frac{1}{\varphi(q)} \; \sideset{}{^{(1)}} \sum_{\xi \bmod{q^2}} \left| \sum_{n=M+1}^{M+N} a_n \xi (n) \right|^2,
\end{equation*}
From \eqref{inieq}, we have
\begin{eqnarray*}
T(N, Q) & = & \sum_q \mathop{\sum_n \sum_{n'}}_{n \equiv n' \bmod{q}} \overline{a}_n a_{n'} e \left( \frac{\overline{n'}(nn'-1)}{q^2} \right) \\
& = & \sum_n \sum_{n'} \sum_{q|(n-n')} \overline{a}_n a_{n'} e \left( \frac{\overline{n'}(nn'-1)}{q^2} \right) \\
& = & \sum_{|l| \leq 2N/Q} \mathop{\sum_q \sum_n \sum_{n'}}_{ql=n-n'} \overline{a}_n a_{n'} e \left( \frac{\overline{n}l}{q} \right).
\end{eqnarray*}
We get, by applying Cauchy's inequality,
\begin{eqnarray*}
T^2(N,Q) & \leq & \left( \sum_{n'} |a_{n'}|^2 \right) \sum_{n'} \left| \mathop{\sum \sum \sum}_{n-n'=lq} a_n e \left( \frac{\overline{n} l}{q} \right) \right|^2  \\
 & = & \left( \sum_{n'} |a_{n'}|^2 \right) \mathop{\sum \sum \sum \sum \sum \sum}_{n_1-n_2=l_1q_1-l_2q_2} a_{n_1} \overline{a}_{n_2} e \left( \frac{\overline{n}_1l_1}{q_1} - \frac{\overline{n}_2l}{q_2} \right).
\end{eqnarray*}
From the above, we get that
\begin{equation*}
T^2(N,Q) \leq \left( \sum_n |a_n|^2 \right) \sum_{n_1} \sum_{n_2} |a_{n_1} a_{n_2} | \sum_{q_1} \left| \mathop{\sum \sum \sum}_{l_1q_1-l_2q_2=n_1-n_2} e \left( \frac{\overline{(n_2-l_2q_2)}l_1}{q_1}-\frac{\overline{n}_2l_2}{q_2} \right) \right|.
\end{equation*}

Here, we apply once again, the Cauchy-Schwartz inequality, getting
\begin{equation*}
T^4(N,Q) \leq \left( \sum_n |a_n|^2 \right)^4 Q \sum_{n_1} \sum_{n_2} \sum_q \left| \mathop{\sum_{q_1} \sum_{l_1} \sum_{l}}_{l_1q_1-l_2q_2=n_1-n_2} e \left( \frac{\overline{(n_1-lq)}l}{q_1}-\frac{\overline{n}_2l}{q} \right) \right|^2.
\end{equation*}
Opening the square modulus, we get
\begin{equation} \label{incompkloo}
T^4(N,Q) \leq \left( \sum_n |a_n|^2 \right)^4 Q \mathop{\sum_q \sum_{q_1} \sum_{q_2} \sum_{l_1} \sum_{l_2} \sum_{l'} \sum_{l''}}_{l_1q_1-l_2q_2=q(l'-l'')} \left\{ \sum_n \right\},
\end{equation}
where the inner-most sum is
\begin{equation*}
\sum_n e \left[ \frac{\overline{(n-l'q)}l_1}{q_1}-\frac{\overline{n}l'}{q}-\frac{\overline{(n-l''q)}l_2}{q_2}+\frac{\overline{n}l''}{q} \right],
\end{equation*}
which is an incomplete Kloosterman type sum and may be completed by Fourier techniques. \newline

It suffices to estimate a sum of the following form
\begin{equation} \label{expoamp}
\sum_{N \leq n \leq N_1} e \left[ \frac{\overline{(n-l'q)}l_1}{q_1}-\frac{\overline{n}l'}{q}-\frac{\overline{(n-l''q)}l_2}{q_2}+\frac{\overline{n}l''}{q} \right],
\end{equation}
where $N < N_1 \leq 2N$.  Our sum may be written as
\begin{equation} \label{incompklooster}
\frac{1}{q_1q_2q} \sum_{a \bmod{q_1q_2q}} \sum_n e \left( \frac{an}{q_1q_2q} \right) \sum_{x \bmod{q_1q_2q}} e \left[ f(x) - \frac{ax}{q_1q_2q} \right],
\end{equation}
where $f(x)$ is the amplitude in \eqref{expoamp} and the range of summation for $n$ is the same as before.  The main contribution will come from the part where $a \equiv 0 \pmod{q_1q_2q}$.  That part in \eqref{incompklooster} is the following
\begin{equation} \label{lsram}
\frac{N}{q_1q_2q}  \sum_{x \bmod{q_1q_2q}} e \left[ f(x) \right].
\end{equation}

The other parts will have small contributions, but they are nevertheless
\begin{equation}
\label{lskloo} \ll \sum_{0<|a| \leq q_1q_2q/2} \frac{1}{|a|} \left| \sum_{x \bmod{q_1q_2q}} e \left[ f(x) - \frac{ax}{q_1q_2q} \right] \right|,
\end{equation}
where we have applied the bound for the geometric series over $n$.  The sum in \eqref{lsram} is a Ramanujan type sum while the one in \eqref{lskloo} may be estimated via Weil's bound for Kloosterman sums, Theorem~\ref{weilbound}. \newline

Making the change of variable from $x$ into $x+l''q$, the inner-most sum of \eqref{lskloo}, it becomes
\begin{equation*}
e \left( \frac{al''}{qq_1q_2} \right) \sum_{x \bmod{q_1q_2q}} e \left[ \left( \frac{\overline{(x+lq)}l_1}{q_1} - \frac{\overline{x}l_2}{q_2} + \frac{\overline{x}l}{q} \right) - \frac{ax}{q_1q_2q} \right],
\end{equation*}
where $l=l''-l'$.  The above is, after a change of variables, a complete Kloosterman sum of modulus $qq_1q_2$.  Hence we apply Weil's bound for Kloosterman sums and the above is $ \ll \gcd(a, q_1q_2q)^{\frac{1}{2}} Q^{\frac{3}{2}+\epsilon}$.  Recall that we are assuming that $q_1$, $q_2$ and $q$ are primes.  Summing over $a$, we get that the sum in \eqref{lskloo} is $\ll Q^{\frac{3}{2}+\epsilon}$. \newline

It still remains to estimate the sum in \eqref{lsram}.  The sum simplifies to
\begin{equation} \label{lsram2}
\sum_{x \bmod{q_1q_2q}} e \left[ \frac{\overline{(x+lq)}l_1}{q_1}-\frac{\overline{x}l_2}{q_2}+ \frac{\overline{x}l}{q} \right].
\end{equation}
If $q_1$, $q_2$ and $q_3$ are pair-wise co-prime, {\it id est} distinct, we write
\begin{equation*}
x = x_1 \overline{q_2q}q_2q+x_2\overline{q_1q}q_1q + y \overline{q_1q_2}q_1q_2,
\end{equation*}
where $x_1 \pmod{q_1}$, $x_2 \pmod{q_2}$ and $y \pmod{q}$.  The sum of our interest in \eqref{lsram2} factors into the product of three Ramanujan type sums
\begin{equation*}
R \left( \frac{l_1}{q_1} \right) R \left( \frac{l_2}{q_2} \right) R \left( \frac{l}{q} \right), \; \mbox{where} \;
R \left( \frac{l}{q} \right) = \sum_{\substack{n \bmod{q} \\ \gcd(n,q)=1}} e \left( \frac{nl}{q} \right), \; \mbox{and} \; \left| R \left( \frac{l}{q} \right) \right| \leq \gcd(l,q)
\end{equation*}
by the virtue of \eqref{ramsumest} and we also note the following almost trivial estimate.
\begin{equation*}
\sum_{d \bmod{q}} \gcd(d,q) \leq \sum_{l|q} \frac{q}{l}l = q \tau(q) \ll q^{1+\epsilon}.
\end{equation*}
Hence, summing over all relevant variables, we have the contribution of \eqref{lsram} to the majorant of $T^4(N,Q)$ is
\begin{equation*}
QN \left( \frac{N}{Q} + Q \right)^3 \left( \sum_n |a_n|^2 \right)^4.
\end{equation*}
If some of the $q$'s are not pair-wise distinct, then the estimate will essentially go the same way as before, but the contribution to the majorant will be different.  If two of the $q$'s are the same, then the contribution of the majorant is
\begin{equation*}
Q^{\epsilon} \left( N^4 Q^{-2} + N^2 Q^2 \right) \left( \sum_n |a_n|^2 \right)^4,
\end{equation*}
and similarly if all $q$'s are the same, the contribution is
\begin{equation*}
Q^{\epsilon} \left( N^4Q^{-2} + N^3 \right) \left( \sum_n |a_n|^2 \right)^4.
\end{equation*}
Combining everything, the left-hand side of \eqref{lsspeceq} is
\begin{equation} \label{semires}
\ll Q^{\epsilon} \left( NQ^{\frac{1}{2}} + N^{\frac{1}{4}} Q^2 + N^{\frac{1}{2}} Q^{\frac{3}{2}} + N^{\frac{3}{4}} Q^{\frac{11}{8}} \right) \sum_n |a_n|^2.
\end{equation}
The last term in the majorant of \eqref{semires} comes from the contribution of \eqref{lskloo}.  The third term in \eqref{semires} is not needed, as if $N^{\frac{1}{2}} Q^{\frac{3}{2}} \gg N^{\frac{1}{4}} Q^2$, then $N^{\frac{1}{4}} \gg Q^{\frac{1}{2}}$.  In that case, $N^{\frac{3}{4}}Q^{\frac{11}{8}} \gg N^{\frac{3}{4}}Q \gg N^{\frac{1}{2}} Q^{\frac{3}{2}}$.  Our contention follows. \end{proof}

\section{Notes}

The asymptotic formula of \eqref{lsspectriv} is useful when $N \ll Q$.  Theorem~\ref{lsspec} is better than \eqref{lsspectriv} when $Q^{\frac{3}{2}+\epsilon} \ll N$.  We would certainly hope to have a result that is useful whenever \eqref{lsspectriv} is not, {\it id est} whenever $Q \ll N$.  However, since we already have to resort to the strength of Weil bound for our present result, any desire for improvement is perhaps too greedy.

\bibliography{biblio}
\bibliographystyle{amsxport}

\vspace*{.1in}
\hspace*{.5in}{\sc\small Dept. Math., U.S. Military Academy, West Point, NY 10996 \newline
\hspace*{.5in}\indent Email Address: {\tt al1526@usma.edu} \newline
\hspace*{.5in}\indent Webpage: {\tt http://www.dean.usma.edu/math/People/Zhao/}}

\end{document}